\let\OToneAccents=\relax
\magnification=\magstep1
\frenchspacing
\baselineskip=16truept
\font\douze=cmr10 at 12pt

\font\grande=cmb10 at 16truept
\def\titre#1{{\OToneAccents\noindent\grande #1}}

\def \chapter#1{\vfill\eject\ifodd\pageno\else\ \vfill\eject\fi\centerline{\grande #1}\bigskip}

\font\tendo=wncyr10 at 12pt
\font\sevendo=wncyr7
\newfam\dofam 
\textfont\dofam = \tendo 
\scriptfont\dofam= \sevendo

\font\tendb=msbm10 
\font\sevendb=msbm7
\newfam\dbfam 
\textfont\dbfam = \tendb 
\scriptfont\dbfam= \sevendb
\def\db {\fam\dbfam\tendb}
\font\tenrsfs=rsfs10 
\font\sevenrsfs=rsfs7
\newfam\scrfam 
\textfont\scrfam = \tenrsfs
\scriptfont\scrfam= \sevenrsfs

\def\Q{{\db Q}}

\def\Z{{\db Z}}

\def\Gal{\mathop{\rm Gal}\nolimits}

\font\tendb=msbm10 
\font\sevendb=msbm7
\newfam\dbfam 
\textfont\dbfam = \tendb 
\scriptfont\dbfam= \sevendb
\def\db {\fam\dbfam\tendb}
\font\tenrsfs=rsfs10 
\font\sevenrsfs=rsfs7
\newfam\scrfam 
\textfont\scrfam = \tenrsfs
\scriptfont\scrfam= \sevenrsfs

\def\Q{{\db Q}}
\def\Z{{\db Z}}

\def\Gal{\mathop{\rm Gal}\nolimits}
\def\mod{\hbox{\rm \ mod.\ }}

\font\teufm=eufm10
\font\seufm=eufm10 at 7pt
\font\sseufm=eufm10 at 6pt
\newfam\fameufm
\textfont\fameufm=\teufm
\scriptfont\fameufm=\seufm
\scriptscriptfont\fameufm=\sseufm
\def\goth{\fam\fameufm\teufm}

\def\qq{{\goth q}}

\def\bb{{\goth b}}

\def\aa{{\goth a}}

\def\Gal{\mathop{\rm Gal}\nolimits}

\def\mod{\hbox{\rm \ mod.\ }}

\centerline {\titre{\'Equation de Fermat et nombres premiers inertes }}
\bigskip
\centerline {Alain Kraus}
\bigskip 
\smallskip
{\bf{Abstract.}}  Let $K$ be a number field and $p$  a prime number $\geq 5$. Let us denote by $\mu_p$ the group of the $p$th roots of unity. We define $p$ to be $K$-regular if 
$p$ does not divide the class number of the field $K(\mu_p)$. 
Under the assumption that $p$ is $K$-regular  and inert in $K$, we establish the  second case of Fermat's Last Theorem over $K$ for the exponent $p$. 
We use in the proof classical arguments, as well as  Faltings' theorem  stating  that a curve of genus at least two over $K$ has a finite number of $K$-rational points.  Moreover, if $K$ is an imaginary quadratic field, other than $\Q(\sqrt{-3}\bigr)$, we deduce a statement which allows often in practice to prove Fermat's Last Theorem 
over $K$ for the  $K$-regular exponents.
\bigskip

\medskip

\bigskip
\medskip

\centerline{\douze{Introduction}}
\bigskip
Soient $K$ un corps de nombres  et $p$ un nombre premier $\geq 5$. On dit que le   th\'eor\`eme de Fermat est vrai sur $K$ pour 
l'exposant $p$, s'il n'existe pas  d'\'el\'ements $x,y,z$ dans $K$ tels que  l'on ait
$$x^p+y^p+z^p=0\quad \hbox{et}\quad  xyz\neq 0.$$

 A. Wiles a d\'emontr\'e  en 1994 qu'il en est ainsi pour le corps $\Q$ ([10]).  Ce r\'esultat a \'et\'e \'etendu au corps $\Q\bigl(\sqrt{2}\bigr)$  
 par 
F. Jarvis et P. Meekin en 2004 ([7]), et tout r\'ecemment aux corps    $\Q\bigl(\sqrt{d}\bigr)$  avec  $d$  sans facteurs carr\'es, $3\leq d\leq 23$ et $d\neq 5,17$, par N. Freitas et S. Siksek ([4]).  Dans le cas o\`u  $K$ est un corps totalement r\'eel, 
ils ont aussi 
d\'emontr\'e  que si  une certaine condition est satisfaite par $K$, le th\'eor\`eme de Fermat est vrai sur $K$ pour tout  $p$  assez grand ([3]). Ce sont les principaux r\'esultats obtenus sur l'\'equation de Fermat au cours de ces vingt derni\`eres ann\'ees.
\smallskip
Notons  $\mu_p$ le groupe des racines $p$-i\`emes de l'unit\'e dans une cl\^oture alg\'ebrique de $K$. Adoptons la terminologie selon laquelle  $p$ est $K$-r\'egulier si 
$p$ ne divise pas le nombre de classes de $K(\mu_p)$. Si  $p$ est $K$-r\'egulier et inerte dans $K$, 
on  d\'emontre ici le  second cas du th\'eor\`eme de Fermat sur $K$ pour l'exposant $p$. On   utilise pour cela   des techniques classiques (cf. [9], chap. 9, si $K=\Q$),  ainsi que  le th\'eor\`eme de G. Faltings sur la finitude de l'ensemble  des points $K$-rationnels  des courbes d\'efinies sur $K$ de genre au moins deux ([2]).

Si  $K$ est un corps quadratique,  F.~H. Hao et C. J. Parry ont obtenu en 1984 des caract\'e\-risations simples des nombres premiers $K$-r\'eguliers ([6]). 
De plus, ils ont \'etabli le  second cas du th\'eor\`eme de Fermat sur $K$ pour les exposants $K$-r\'eguliers et d\'ecompos\'es dans $K$ ([5], th. 3).
Si $K$ est un corps quadratique imaginaire, autre que $\Q\bigl(\sqrt{-3}\bigr)$, 
avec les r\'esultats obtenus  r\'ecemment  dans  [8],  
on en d\'eduit un \'enonc\'e permettant souvent en pratique de d\'emontrer le th\'eor\`eme de Fermat sur $K$  pour les exposants $K$-r\'eguliers.
\smallskip
Je remercie  D.  Bernardi pour les conversations que nous avons eues concernant cet article.
\bigskip

{\bf{ 1. \'Enonc\'e des r\'esultats}}
\medskip
Soient $K$ un corps de nombres, d'anneau d'entiers $O_K$, et $p$ un nombre premier $\geq 5$.  
\bigskip

\proclaim Th\'eor\`eme.  Supposons que  $p$ soit $K$-r\'egulier et inerte dans $K$. 
Alors,  il n'existe pas de triplets $(x,y,z)$ d'\'el\'ements de $O_K$ tels que 
$$x^p+y^p+z^p=0,\quad xyz\neq 0,\quad (xy)O_K+pO_K=O_K\quad \hbox{et}\quad p\ \hbox{divise}\ z.\leqno(1)$$
\medskip

Pour tout $n\geq 1$, notons  $W_n$ le r\'esultant des polyn\^omes  $X^n-1$ et $(X+1)^n-1$. 
\bigskip

\proclaim Proposition. Supposons  que $K$ soit un corps quadratique imaginaire et les deux conditions suivantes satisfaites :
\smallskip
\vskip0pt\noindent
1)  il existe  $n\geq 1$ tel que  $np+1$ soit un nombre premier d\'ecompos\'e dans $K$  ne divisant pas $(n^n-1)W_n$.
\smallskip
\vskip0pt\noindent
2) $p$ est   $K$-r\'egulier et  non ramifi\'e dans $K$.
\smallskip
\vskip0pt\noindent
Alors, le   th\'eor\`eme de Fermat est vrai sur $K$ pour l'exposant $p$.
\bigskip

On en d\'eduit par exemple que le th\'eor\`eme de Fermat est vrai sur  le corps $\Q(i)$ pour les nombres premiers  $p<10^6$ qui sont $\Q(i)$-r\'eguliers ([8], cor. 2).  La liste de ceux plus petits que  100 est 
$$5,7,11,13,17,23,29,41,53,73,83,89,97.$$
Il semble donc que 19 soit le plus petit nombre premier   pour lequel on ne sache pas conclure sur $\Q(i)$.

\bigskip
{\bf{ 2.  Notations}}
\medskip
Pour tout corps de nombres $N$, on notera $O_N$ son anneau d'entiers. 
 Soit $\zeta$ un g\'en\'e\-rateur de $\mu_p$.  Posons 
$$L=K(\mu_p),\quad \pi=1-\zeta\quad \hbox{et}\quad \lambda=(1-\zeta)(1-\zeta^{-1}).$$

Le nombre premier $p$ \'etant non ramifi\'e dans $K$ et totalement ramifi\'e dans $\Q(\mu_p)$, on a $K\cap \Q(\mu_p)=\Q$. On a donc les \'egalit\'es des degr\'es
$$[L:K]=p-1\quad \hbox{et}\quad [L:K(\lambda)]=2.$$
Par hypoth\`ese,  $pO_K$ est un id\'eal premier de $O_K$. Il est totalement ramifi\'e dans $L$ et  $\pi O_L$ est donc l'unique id\'eal premier de $O_L$ au-dessus de $pO_K$.  De m\^eme  
$\lambda O_{K(\lambda)}$ est l'unique id\'eal premier de $O_{K(\lambda)}$ au-dessus de $pO_K$, et on a l'\'egalit\'e
$$\lambda O_L=(\pi O_L)^2.\leqno(2)$$
On notera $v_{\lambda}$ (resp. $v_{\pi}$) la valuation de $K(\lambda)$ (resp.  $L$) associ\'ee \`a $\lambda$ (resp. $\pi$).

\bigskip
{\bf{ 3.  Principe de d\'emonstration du th\'eor\`eme}}
\medskip
Soit $S$ l'ensemble  des quintuplets $(\alpha,\beta,\gamma,\mu,n)$ v\'erifiant les conditions suivantes :
\smallskip
1) $\alpha,\beta,\gamma$ sont des \'el\'ements de  $O_{K(\lambda)}$ non divisibles par $\lambda$.
\smallskip
2)  $\mu$ est une unit\'e de $O_{K(\lambda)}$ et $n$ est un entier $\geq 2$.
\smallskip
3) On a l'\'egalit\'e
 $$\alpha^p+\beta^p=\mu \bigl(\lambda^n\gamma\bigr)^p.\leqno(3)$$
\medskip 
On proc\`ede par l'absurde en supposant qu'il existe un triplet d'\'el\'ements de $O_K$ satisfaisant la condition (1).  
Cela entra\^\i ne  que  $S$ n'est pas vide. Soit $(\alpha,\beta,\gamma,\mu,n)$ un \'el\'ement de   $S$.  On d\'emontre alors  qu'il existe  $\alpha',\beta',\gamma'$ dans $O_{K(\lambda)}$, non divisibles par $\lambda$, et  une unit\'e $\mu'$ de $O_{K(\lambda)}$ tels que 
$$(\alpha',\beta',\gamma',\mu',2n-1)$$
soit aussi un \'el\'ement de $S$.
On conclut  comme suit. 
 Parce que $n\geq 2$, on a 
 $$2n-1>n.\leqno(4)$$
Soit   $U$ un syst\`eme de repr\'esentants des unit\'es de $O_{K(\lambda)}$ modulo les puissances $p$-i\`emes. D'apr\`es  la condition (4), l'ensemble des quadruplets
   $(x,y,z,u)$ tels que  l'on ait
   $$x^p+y^p=uz^p\quad \hbox{avec}\quad x,y,z\in K(\lambda)\quad \hbox{et}\quad u\in U$$
   est infini.  L'ensemble  $U$ est fini, donc il existe   $u\in U$ tel que la courbe  d'\'equation
 $$x^p+y^p=uz^p$$
 poss\`ede une infinit\'e de points rationnels sur $K(\lambda)$. Son genre \'etant au moins deux, cela contredit  le th\'eor\`eme de Faltings. Ainsi, $S$ est vide, d'o\`u une contradiction et le th\'eor\`eme annonc\'e.
\smallskip
 L'essentiel de la suite est consacr\'e \`a  sa d\'emonstration.
\bigskip
 {\bf{4. Lemmes pr\'eliminaires}}
 \medskip
Supposons d\'esormais  qu'il existe $x,y,z$ dans $O_K$ v\'erifiant la condition (1).
\smallskip

\proclaim Lemme 1. L'ensemble $S$ n'est pas vide.
 
D\'emonstration : Il existe  $z_1$ dans $O_K$, non divisible par $p$, 
et un entier $a\geq 1$  tels que  
$$x^p+y^p=-z^p=(p^az_1)^p.$$
 Par ailleurs, il existe une unit\'e $\varepsilon\in \Z[\lambda]$ telle que l'on ait
 $$p=\varepsilon \lambda^{{p-1\over 2}}.$$
En posant $\eta=\varepsilon^{ap}$, on obtient 
 $$x^p+y^p=\eta\Bigl(\lambda^{{a(p-1)\over 2}} z_1\Bigr)^p.$$
 On a $p\geq 5$, d'o\`u 
 $${a(p-1)\over 2}\geq {p-1\over 2}\geq 2.$$
 Par suite, le quintuplet $\Bigl(x,y,z_1,\eta,{a(p-1)\over 2}\Bigr)$ est dans $S$, d'o\`u le lemme.
\bigskip

Consid\'erons, pour toute la suite,    un \'el\'ement $(\alpha,\beta,\gamma,\mu,n)$ de  $S$. D'apr\`es (3),  on a 
 $$\prod_{i=0}^{p-1}(\alpha+\beta \zeta^i)=\mu (\lambda^n\gamma)^p.\leqno(5)$$

\proclaim Lemme 2. Pour tout $i=0,\cdots,p-1$, on a  $\alpha+\beta \zeta^i\equiv 0 \mod \pi O_L$.

D\'emonstration :  D'apr\`es (5), $\pi$ divise l'un des $\alpha+\beta \zeta^i$ et  pour tout $j$ on a 
$$\alpha+\beta \zeta^j\equiv \alpha+\beta \mod \pi O_L,$$
d'o\`u l'assertion.
\bigskip

Posons
 $$D=\alpha O_L+\beta O_L$$
 et pour tous $i$ et $j$ distincts entre $0$ et $p-1$,
 $$D_{i,j}=\biggl({\alpha+\beta \zeta^i\over \pi}\biggr)O_L+\biggl({\alpha+\beta \zeta^j\over \pi}\biggr)O_L.$$
\medskip

\proclaim Lemme 3. On a $D=D_{i,j}$.

D\'emonstration :   Soient $\qq$ un id\'eal premier non nul de $O_L$ et $r$ un entier $\geq 1$. Supposons  que $\qq^r$ divise $D_{i,j}$. 
Alors, $\qq^r$ divise
${\beta(\zeta^i-\zeta^j)\over \pi}$. On a $i\neq j$, donc 
$\qq^r$ divise $\beta$. De m\^eme, $\qq^r$ divise $\alpha$, ainsi $D_{i,j}$ divise $D$.
Inversement, si  $\qq^r$ divise $D$, alors $\qq$ est distinct de $\pi O_L$ car $\alpha\beta$ est premier avec $\pi$.  L'\'egalit\'e
$$(\pi O_L)D_{i,j}=(\alpha+\beta \zeta^i)O_L+(\alpha+\beta \zeta^j)O_L$$
entra\^\i ne alors  que $\qq^r$ divise $D_{i,j}$, donc $D$ divise $D_{i,j}$.
\bigskip

\proclaim Lemme 4.  Il existe des id\'eaux $\bb_0,\cdots,\bb_{p-1}$ de $O_L$, non divisibles par $\pi O_L$, premiers entre eux deux \`a deux, tels que l'on ait
$$(\alpha+\beta)O_L=(\lambda O_L)^{np-{p-1\over 2}} D \bb_0^p,\leqno(6)$$
$$\biggl({\alpha+\beta \zeta^i\over \pi}\biggr)O_L=D\bb_i^p \quad \hbox{pour tout}\quad i=1,\cdots,p-1.\leqno(7)$$

D\'emonstration : D'apr\`es le lemme 2,  
$\alpha+\beta$ appartient \`a  $\pi O_L\cap O_{K(\lambda)}$ qui est  $\lambda O_{K(\lambda)}.$
Ainsi, l'\'egalit\'e (2) implique  
$$\alpha+\beta\equiv 0 \mod \pi^2O_L.$$
Pour tout $i$, on a donc  la congruence
$$\alpha+\beta \zeta^i\equiv \alpha(1-\zeta^i) \mod \pi^2O_L.$$
On a $v_{\pi}(\alpha)=0$,  d'o\`u 
pour tout $i=1,\cdots,p-1$, 
$$v_{\pi}(\alpha+\beta \zeta^i)=1.$$
D'apr\`es (5), on obtient
$$v_{\pi}(\alpha+\beta)=2np-(p-1),$$
autrement dit,
$$v_{\lambda}(\alpha+\beta)=np-{p-1\over 2}.$$
Il  existe donc un id\'eal $\aa_0$ de $O_L$, non divisible par $\pi$, tel que 
$$(\alpha+\beta)O_L=(\lambda O_L)^{np-{p-1\over 2}} D \aa_0.$$
De m\^eme, pour tout $i=1,\cdots,p-1$, il existe un id\'eal $\aa_i$ de $O_L$, non divisible par $\pi$,  tel que 
$$(\alpha+\beta \zeta^i)O_L=(\pi O_L) D \aa_i.$$
On a ainsi l'\'egalit\'e
$$(\lambda O_L)^{np} (\gamma O_L)^p=\prod_{i=0}^{p-1}(\alpha+\beta \zeta^i)O_L=(\pi O_L)^{p-1} (\lambda O_L)^{np-{p-1\over 2}}   D^p\aa_0\cdots\aa_{p-1},$$
d'o\`u 
$$(\gamma O_L)^p=D^p\aa_0\cdots\aa_{p-1}.$$
D'apr\`es le lemme 3, les id\'eaux $\aa_i$ sont premiers entre eux deux \`a deux. Chaque id\'eal $\aa_i$ est donc  la puissance $p$-i\`eme d'un id\'eal $\bb_i$ de $O_L$, d'o\`u le lemme.
\bigskip

{\bf{5. Les \'el\'ements $\Theta_k$ de $L$}}
\medskip
Pour tout entier $k$ non multiple de $p$, posons 
$$\Theta_k=\Biggl({\alpha+\beta \zeta^k\over 1-\zeta^k}\biggr)\Biggl({\alpha+\beta \zeta^{-k}\over 1-\zeta^{-k}}\biggr)^{-1}\in L^*.$$
 On a 
$$\biggl({\alpha+\beta \zeta^k\over \pi}\biggr)O_L=\biggl({\alpha+\beta \zeta^k\over 1-\zeta^k}\biggr)O_L.$$
D'apr\`es la formule (7), en r\'eindexant les id\'eaux $\bb_k$ avec $k$ entre $-{p-1\over 2}$ et ${p-1\over 2}$, 
on a 
$$\Theta_k O_L=\bigl(\bb_k\bb_{-k}^{-1}\bigr)^p.\leqno(8)$$ 
\bigskip

\proclaim Lemme 5. L'\'el\'ement $\Theta_k$ est une puissance $p$-i\`eme dans $L$.

D\'emonstration : V\'erifions  que l'extension $L\bigl({\root p\of {\Theta_k}}\bigr)/L$ est partout non ramifi\'ee, ce qui permet de conclure avec l'hypoth\`ese que $p$ ne divise pas le nombre de classes de $L$. 
On  a 
$$\Theta_k-1=-\zeta^{-k}{(1+\zeta^k)(\alpha+\beta)\over  \alpha(1-\zeta^{-k})+\zeta^{-k}(\alpha+\beta)} .$$
D'apr\`es les \'egalit\'es (2) et (6), on a ainsi
$$v_{\pi}(\Theta_k-1)=v_{\pi}(\alpha+\beta)-1=(2n-1)p.$$
En particulier, on a $\Theta_k\equiv 1 \mod \pi^p$. Il en r\'esulte que l'extension $L\bigl({\root p\of {\Theta_k}}\bigr)/L$ est non ramifi\'ee en  $\pi O_L$ (cf.  [1], cor. p. 503 ; avec ses notations, on a $e(\pi O_L/p)=p-1$ et $z(\pi O_L,p)=p+1$).
Par ailleurs, si $\qq$ est un id\'eal premier de $O_L$ distinct de $\pi O_L$, on a $v_{\qq}(pO_L)=0$, et d'apr\`es la formule (8) on a 
$$v_{\qq}(\Theta_k)\equiv 0 \mod p.$$
Cela entra\^\i ne que $\qq$ est non ramifi\'e dans $L\bigl({\root p\of {\Theta_k}}\bigr)$  ({\it{loc. cit.}}), d'o\`u le r\'esultat.
\vfill\eject

{\bf{6. Les unit\'es $\varepsilon_k$ de $O_{K(\lambda)}$ et les entiers $\rho_k$ de $O_L$ }}
\medskip 
Notons $h_{\lambda}$ le nombre de classes de $K(\lambda)$. Par hypoth\`ese,  $p$ est $K$-r\'egulier. Le degr\'e de $L$ sur $K(\lambda)$ est premier \`a $p$,  donc $p$ ne divise pas $h_{\lambda}$. Il existe ainsi  un entier $t\geq 1$ tel que 
$$th_{\lambda} \equiv -2 \mod p.\leqno(9)$$
L'id\'eal $\bigl(\alpha O_{K(\lambda)}+\beta O_{K(\lambda)}\bigr)^{th_{\lambda}}$ \'etant principal, il existe  $d\in O_{K(\lambda)}$ tel que l'on ait
$$\bigl(\alpha O_{K(\lambda)}+\beta O_{K(\lambda)}\bigr)^{th_{\lambda}}=dO_{K(\lambda)}.$$
On a ainsi l'\'egalit\'e
$$D^{th_{\lambda}}=dO_L.\leqno(10)$$
D\'esignons dans la suite par $\sigma$ l'\'el\'ement non trivial du groupe de Galois $\Gal\bigl(L/K(\lambda)\bigr)$. On~a 
$$\sigma(\zeta)=\zeta^{-1}.$$
\medskip

\proclaim Lemme 6. Soit $k$ un entier non multiple de $p$. Il existe une unit\'e $\varepsilon_k$ de $O_{K(\lambda)}$ et un \'el\'ement $\rho_k$ de $O_L$ tels que l'on ait 
$$d^{{p+1\over 2}} \biggl({\alpha+\beta \zeta^k\over 1-\zeta^k}\biggr)=\varepsilon_k^{{p+1\over 2}}\rho_k^p.\leqno(11)$$

D\'emonstration :  Utilisons  la formule (7) avec $i=k$ et $i=-k$.  On obtient   
$$\biggl({\alpha+\beta \zeta^k\over 1-\zeta^k}.{\alpha+\beta \zeta^{-k}\over 1-\zeta^{-k}}\biggr)O_L=D^2\bigl(\bb_k\bb_{-k}\bigr)^p.$$
D'apr\`es (9) et (10), on a donc
$$(dO_L)\biggl({\alpha+\beta \zeta^k\over 1-\zeta^k}.{\alpha+\beta \zeta^{-k}\over 1-\zeta^{-k}}\biggr)O_L=\Bigl(D^{{th_{\lambda}+2\over p}}\bb_k\bb_{-k}\Bigr)^p.\leqno(12)$$
L'\'el\'ement
$$d\biggl({\alpha+\beta \zeta^k\over 1-\zeta^k}.{\alpha+\beta \zeta^{-k}\over 1-\zeta^{-k}}\biggr)$$
est fix\'e par $\sigma$. Il appartient donc \`a $O_{K(\lambda)}$.
On  d\'eduit alors de (12) que l'id\'eal 
$$d\biggl({\alpha+\beta \zeta^k\over 1-\zeta^k}.{\alpha+\beta \zeta^{-k}\over 1-\zeta^{-k}}\biggr) O_{K(\lambda)}$$
est la puissance $p$-i\`eme d'un id\'eal de $O_{K(\lambda)}$.
Puisque $p$ ne divise pas $h_{\lambda}$, il existe donc $\varphi_k\in O_{K(\lambda)}$ et une unit\'e $\varepsilon_k\in O_{K(\lambda)}$ tels que l'on ait
$$d\biggl({\alpha+\beta \zeta^k\over 1-\zeta^k}\biggr)\biggl({\alpha+\beta \zeta^{-k}\over 1-\zeta^{-k}}\biggr)=\varepsilon_k \varphi_k^p.$$
D'apr\`es le lemme 5, il existe  $\Psi_k\in L$ tel que l'on ait
$$\Theta_k=\Psi_k^p.$$
On en d\'eduit l'\'egalit\'e 
$$d\biggl({\alpha+\beta \zeta^k\over 1-\zeta^k}\biggr)^2=\varepsilon_k(\varphi_k \Psi_k)^p.$$
En \'elevant ses deux membres  \`a la puissance ${p+1\over 2}$, on a
$$d^{{p+1\over 2}}\biggl({\alpha+\beta \zeta^k\over 1-\zeta^k}\biggr)^{p+1}=\varepsilon_k^{{p+1\over 2}}\Bigl((\varphi_k \Psi_k)^{{p+1\over 2}}\Bigr)^p.$$
Posons
$$\rho_k=(\varphi_k \Psi_k)^{{p+1\over 2}}\biggl({\alpha+\beta \zeta^{k}\over 1-\zeta^k}\biggr)^{-1}.$$
 C'est un \'el\'ement de $L$. Parce que $\varepsilon_k$ est une unit\'e de $O_{K(\lambda)}$, $\rho_k^p$ est dans $O_L$, donc $\rho_k$ aussi, d'o\`u l'\'egalit\'e (11).
 \bigskip

\proclaim Lemme 7. Il existe une unit\'e $\varepsilon_0$ de $O_{K(\lambda)}$ et un \'el\'ement $\rho_0$ de $O_{K(\lambda)}$ tels que l'on ait
$$d(\alpha+\beta)^2=\varepsilon_0 \lambda^{2np-(p-1)}\rho_0^p.$$

D\'emonstration :   D'apr\`es la formule (6), on a 
$$(\alpha+\beta)^2O_L= \bigl(\lambda O_L\bigr)^{2np-(p-1)}D^2\bb_0^{2p}.$$
En utilisant (9) et (10), on a donc
$$d(\alpha+\beta)^2O_L=\bigl(\lambda O_L\bigr)^{2np-(p-1)}    \Bigl(D^{{th_{\lambda}+2\over p}}\bb_0^2\Bigr)^p.$$
L'\'el\'ement 
$${d(\alpha+\beta)^2\over \lambda^{2np-(p-1)}}$$
\'etant dans $O_{K(\lambda)}$, l'id\'eal 
$$\biggl({d(\alpha+\beta)^2\over \lambda^{2np-(p-1)}}\biggr)O_{K(\lambda)}$$
est donc la puissance $p$-i\`eme d'un id\'eal de $O_{K(\lambda)}$.
Le fait que  $p$ ne divise pas $h_{\lambda}$ implique alors le r\'esultat.
\bigskip

{\bf{7. Fin de la d\'emonstration}}
\medskip 
Soit $k$ un entier non multiple de $p$. Parce que l'on a $p\geq 5$, il existe un entier $\ell$ tel que l'on ait $\ell\not\equiv 0,\pm k \mod p$.
\bigskip

\proclaim Lemme 8. Il existe  une unit\'e $\delta$ de $O_{K(\lambda)}$ telle  que l'on ait
$$\biggl({\varepsilon_k\over \varepsilon_{\ell}}\biggr) \bigl(\varepsilon_k\rho_k\sigma(\rho_k)\bigr)^p+\bigl(-\varepsilon_{\ell}\rho_{\ell}\sigma(\rho_{\ell})\bigr)^p=\delta \bigl(\lambda^{2n-1}\rho_0d\bigr)^p.\leqno(13)$$

D\'emonstration :  En composant   l'\'egalit\'e (11) par $\sigma$,
on a 
$$d^{{p+1\over 2}} \biggl({\alpha+\beta \zeta^{-k}\over 1-\zeta^{-k}}\biggr)=\varepsilon_k^{{p+1\over 2}} \sigma(\rho_k)^p.\leqno(14)$$
Pour tout $i$,  posons dans $O_{K(\lambda)}$
$$\lambda_i=(1-\zeta^i)(1-\zeta^{-i}).$$
En effectuant le produit des \'egalit\'es (11) et (14), on obtient
$$d^{p+1}  \Bigl((\alpha+\beta)^2-\alpha\beta \lambda_k\Bigr)=\lambda_k \varepsilon_k^{p+1}\bigl(\rho_k\sigma(\rho_k)\bigr)^p.$$
On en d\'eduit avec le lemme 7 l'\'egalit\'e
$$\varepsilon_0 \lambda^{2np-(p-1)}(\rho_0d)^p-d^{p+1}\alpha\beta \lambda_k = \lambda_k \varepsilon_k^{p+1}\bigl(\rho_k\sigma(\rho_k)\bigr)^p.\leqno(15)$$
On a de m\^eme l'\'egalit\'e
$$\varepsilon_0 \lambda^{2np-(p-1)}(\rho_0d)^p-d^{p+1}\alpha\beta \lambda_{\ell} = \lambda_{\ell} \varepsilon_{\ell}^{p+1}\bigl(\rho_{\ell}\sigma(\rho_{\ell})\bigr)^p.\leqno(16)$$
Multiplions (15) par $\lambda_{\ell}$ et (16) par $\lambda_k$, puis soustrayons  les deux \'egalit\'es obtenues. On obtient
$$\biggl({\lambda_{\ell}-\lambda_k\over \lambda_{\ell}\lambda_k}\biggr)\bigl(\varepsilon_0 \lambda^{2np-(p-1)}(\rho_0d)^p\bigr)=\varepsilon_k^{p+1}\bigl(\rho_k\sigma(\rho_k)\bigr)^p-\varepsilon_{\ell}^{p+1}\bigl(\rho_{\ell}\sigma(\rho_{\ell})\bigr)^p.$$
Par ailleurs, on a 
$$\lambda_{\ell}-\lambda_k=\bigl(1-\zeta^{k+\ell}\bigr)\bigl(\zeta^{-k}-\zeta^{-\ell}\bigr).$$
On a $\ell\not\equiv \pm k \mod p$, donc 
$1-\zeta^{k+\ell}$ et $\zeta^{-k}-\zeta^{-\ell}$
sont associ\'es \`a $\pi$ dans $O_L$. Ainsi  ${\lambda_{\ell}-\lambda_k\over \lambda_{\ell}\lambda_k}$ et ${1\over \lambda}$ sont  associ\'es   dans $O_L$. 
Ce sont des \'el\'ements de  $O_{K(\lambda)}$, donc il existe une unit\'e $\delta'$ de $O_{K(\lambda)}$ telle que l'on ait
$${\lambda_{\ell}-\lambda_k\over \lambda_{\ell}\lambda_k}={\delta'\over \lambda}.$$
On  en d\'eduit  l'\'egalit\'e
$$\delta'\varepsilon_0 \bigl(\lambda^{2n-1} \rho_0d\bigr)^p=\varepsilon_k\bigl(\varepsilon_k\rho_k\sigma(\rho_k)\bigr)^p+\varepsilon_{\ell}\bigl(-\varepsilon_{\ell}\rho_{\ell}\sigma(\rho_{\ell})\bigr)^p.$$
En posant
$$\delta={\delta'\varepsilon_0\over \varepsilon_{\ell}},$$
qui est  une unit\'e de $O_{K(\lambda)}$, on obtient l'\'egalit\'e (13).
\bigskip

\proclaim Lemme 9. L'unit\'e ${\varepsilon_k\over \varepsilon_{\ell}}$ est  une puissance $p$-i\`eme dans $O_{K(\lambda)}$.

 D\'emonstration :  On a 
$${\alpha+\beta \zeta^k\over 1-\zeta^k}=\alpha+\zeta^k {\alpha+\beta\over 1-\zeta^k}.$$
D'apr\`es la formule (6), on a
$$v_{\pi}\biggl({\alpha+\beta\over 1-\zeta^k}\biggr)=2np-(p-1)-1=(2n-1)p,$$
d'o\`u en particulier la congruence
$${\alpha+\beta \zeta^k\over 1-\zeta^k}\equiv \alpha \mod \pi^p.$$
On d\'eduit alors du lemme 6 que l'on a 
$$d^{p+1} \alpha^2\equiv \varepsilon_k\bigl( \varepsilon_k\rho_k^2\bigr)^p  \mod \pi^p.$$
De m\^eme, on a 
$$d^{p+1} \alpha^2\equiv \varepsilon_{\ell}\bigl( \varepsilon_{\ell}\rho_{\ell}^2\bigr)^p  \mod \pi^p.$$
On a $v_{\pi}(d\alpha)=0$, donc
$\rho_k$ et $\rho_{\ell}$ sont inversibles modulo $\pi^p$.
Par suite, il existe $r\in O_L$ tel que l'on ait
$${\varepsilon_k\over \varepsilon_{\ell}}\equiv r^p\mod \pi^p$$
Il en r\'esulte que l'extension
$$L\biggl( {\root p\of {\varepsilon_k\over \varepsilon_{\ell}}} \biggr)\big/L$$
est partout non ramifi\'ee ([1], cor. p. 503). Le nombre de classes de $L$ n'\'etant pas divisible par $p$,  
${\varepsilon_k\over \varepsilon_{\ell}}$ est  donc une puissance $p$-i\`eme dans $L$, d'o\`u le r\'esultat car ${\varepsilon_k\over \varepsilon_{\ell}}$ est dans $O_{K(\lambda)}$.
\bigskip

Posons
$${\varepsilon_k\over \varepsilon_{\ell}}=\nu^p,$$
o\`u $\nu$ est une unit\'e de $O_{K(\lambda)}$. D'apr\`es le lemme 8, on obtient l'\'egalit\'e
$$ \bigl(\nu\varepsilon_k\rho_k\sigma(\rho_k)\bigr)^p+\bigl(-\varepsilon_{\ell}\rho_{\ell}\sigma(\rho_{\ell})\bigr)^p=\delta \bigl(\lambda^{2n-1}\rho_0d\bigr)^p.$$
 Les \'el\'ements
$$\alpha'=\nu\varepsilon_k\rho_k\sigma(\rho_k),\quad \beta'=-\varepsilon_{\ell}\rho_{\ell}\sigma(\rho_{\ell})\quad \hbox{et}\quad \gamma'=\rho_0d$$
sont dans $O_{K(\lambda)}$, car ils sont fix\'es par $\sigma$. Ils ne sont pas divisibles par $\lambda$ (cf. (6) et (7)). Le quintuplet $(\alpha',\beta',\gamma',\delta,2n-1)$ est donc dans $S$. Compte tenu du paragraphe 3, cela termine la d\'emonstration du th\'eor\`eme.
\bigskip

{\bf{8. La proposition}}
\medskip
Le degr\'e de $K(\mu_p)$ sur $K$ est premier avec $p$. Parce que $p$ est  $K$-r\'egulier, il  ne divise pas le nombre de classes de $K$. La  premi\`ere condition entra\^\i ne ainsi  que le premier cas du th\'eor\`eme de Fermat est vrai sur $K$ pour $p$ ([8]). Le r\'esultat obtenu ici et le th\'eor\`eme 3  de [5] impliquent alors la proposition.

\bigskip
\bigskip  
\centerline {\douze{Bibliographie}}
\bigskip
\vskip0pt\noindent
[1] H. Cohen, Advanced Topics in Computational Algebraic Number Theory, Springer-Verlag  GTM {\bf{193}}, 2000.
\smallskip
\vskip0pt\noindent
[2] G. Faltings, EndlichkeitssŠtze f\"ur abelsche VarietŠten \"uber Zahlk\"orpern, {\it{Invent. Math.}}  {\bf{73}} (1983), 349-366.
\smallskip
\vskip0pt\noindent
[3]  N. Freitas et S. Siksek, The asymptotic Fermat's Last Theorem for five-sixths of real quadratic fields,
arXiv : 1307.3162v3  (16 Jul 2014), 21 pages, \`a para\^\i tre dans la revue Compositio Mathematica.
\smallskip
\vskip0pt\noindent
[4]  N. Freitas et S. Siksek,   Fermat's Last Theorem over some small real quadratic fields, 
arXiv : 1407.4435v1  (16 Jul 2014), 15 pages.
\smallskip
\vskip0pt\noindent
[5] F. H. Hao, C.  J. Parry, The Fermat equation over quadratic fields, {\it J. Number Theory} {\bf 19} (1984),  115-130.
\smallskip
\vskip0pt\noindent
[6]  F. H. Hao, C.  J. Parry, Generalized Bernoulli Numbers and $m$-Regular Primes, {\it{Math. Comp.}} {\bf{43}} (1984), 273-288.
\smallskip
\vskip0pt\noindent
[7] F. Jarvis et P. Meekin, The Fermat equation over $\Q\bigl(\sqrt{2}\bigr)$, {\it J. Number Theory} {\bf 109} (2004),  182-196.
\smallskip
\vskip0pt\noindent
[8] A. Kraus, Remarques sur le premier cas du th\'eor\`eme de Fermat sur les corps de nombres, arXiv:1410.0546v1 (2 octobre 2014), 9 pages,  \`a para\^\i tre dans la revue Acta Arithmetica. 
\smallskip
\vskip0pt\noindent
[9] L. C. Washington, Introduction to Cyclotomic Fields, Springer-Verlag  GTM {\bf{83}}, 1982.
\smallskip
\vskip0pt\noindent
[10]  A. Wiles, Modular elliptic curves and
Fermat's Last Theorem, {\it {Ann. of Math.}} {\bf{141}} (1995), 443-551.
\bigskip

\line {\hfill{27 novembre 2014}}

\item{} Alain Kraus
\item{}Universit\'e de Paris VI, 
\item{} Institut de Math\'ematiques, 
\item{} 4 Place Jussieu, 75005 Paris,  
\item{} France
\medskip
\vskip0pt\noindent
\item{}e-mail : alain.kraus@imj-prg.fr
\smallskip

\bigskip

\bye